\documentclass[11pt]{article}
\usepackage{amsfonts}
\usepackage{amsfonts,mathrsfs,amssymb,amsthm,mathptm}
\usepackage{amsmath,amscd}
\usepackage{mathptm,pslatex}
\usepackage{color}
\usepackage[all]{xy}
\usepackage{graphicx}
\usepackage{fancyhdr}
\usepackage{hyperref}
\usepackage{appendix}

\begin{document}
\newtheorem{Def}{Definition}[section]
\newtheorem{Bsp}[Def]{Example}
\newtheorem{Prop}[Def]{Proposition}
\newtheorem{Theo}[Def]{Theorem}
\newtheorem{Lem}[Def]{Lemma}
\newtheorem{Koro}[Def]{Corollary}
\theoremstyle{definition}
\newtheorem{Rem}[Def]{Remark}

\newcommand{\add}{{\rm add}}
\newcommand{\con}{{\rm con}}
\newcommand{\gd}{{\rm gldim}}
\newcommand{\dm}{{\rm domdim}}
\newcommand{\cdm}{{\rm codomdim}}
\newcommand{\tdim}{{\rm dim}}
\newcommand{\E}{{\rm E}}
\newcommand{\Mor}{{\rm Morph}}
\newcommand{\End}{{\rm End}}
\newcommand{\ind}{{\rm ind}}
\newcommand{\rsd}{{\rm resdim}}
\newcommand{\rd} {{\rm repdim}}
\newcommand{\ol}{\overline}
\newcommand{\overpr}{$\hfill\square$}
\newcommand{\rad}{{\rm rad}}
\newcommand{\soc}{{\rm soc}}
\renewcommand{\top}{{\rm top}}
\newcommand{\pd}{{\rm pdim}}
\newcommand{\id}{{\rm idim}}
\newcommand{\fld}{{\rm fdim}}
\newcommand{\Fac}{{\rm Fac}}
\newcommand{\Gen}{{\rm Gen}}
\newcommand{\fd} {{\rm findim}}
\newcommand{\Fd} {{\rm Findim}}
\newcommand{\Pf}[1]{{\mathscr P}^{<\infty}(#1)}
\newcommand{\DTr}{{\rm DTr}}
\newcommand{\cpx}[1]{#1^{\bullet}}
\newcommand{\D}[1]{{\mathscr D}(#1)}
\newcommand{\Dz}[1]{{\mathscr D}^+(#1)}
\newcommand{\Df}[1]{{\mathscr D}^-(#1)}
\newcommand{\Db}[1]{{\mathscr D}^b(#1)}
\newcommand{\C}[1]{{\mathscr C}(#1)}
\newcommand{\Cz}[1]{{\mathscr C}^+(#1)}
\newcommand{\Cf}[1]{{\mathscr C}^-(#1)}
\newcommand{\Cb}[1]{{\mathscr C}^b(#1)}
\newcommand{\Dc}[1]{{\mathscr D}^c(#1)}
\newcommand{\K}[1]{{\mathscr K}(#1)}
\newcommand{\Kz}[1]{{\mathscr K}^+(#1)}
\newcommand{\Kf}[1]{{\mathscr  K}^-(#1)}
\newcommand{\Kb}[1]{{\mathscr K}^b(#1)}
\newcommand{\DF}[1]{{\mathscr D}_F(#1)}

\newcommand{\Kac}[1]{{\mathscr K}_{\rm ac}(#1)}
\newcommand{\Keac}[1]{{\mathscr K}_{\mbox{\rm e-ac}}(#1)}

\newcommand{\modcat}{\ensuremath{\mbox{{\rm -mod}}}}
\newcommand{\Modcat}{\ensuremath{\mbox{{\rm -Mod}}}}
\newcommand{\Spec}{{\rm Spec}}

\newcommand{\stmc}[1]{#1\mbox{{\rm -{\underline{mod}}}}}
\newcommand{\Stmc}[1]{#1\mbox{{\rm -{\underline{Mod}}}}}
\newcommand{\prj}[1]{#1\mbox{{\rm -proj}}}
\newcommand{\inj}[1]{#1\mbox{{\rm -inj}}}
\newcommand{\Prj}[1]{#1\mbox{{\rm -Proj}}}
\newcommand{\Inj}[1]{#1\mbox{{\rm -Inj}}}
\newcommand{\PI}[1]{#1\mbox{{\rm -Prinj}}}
\newcommand{\GP}[1]{#1\mbox{{\rm -GProj}}}
\newcommand{\GI}[1]{#1\mbox{{\rm -GInj}}}
\newcommand{\gp}[1]{#1\mbox{{\rm -Gproj}}}
\newcommand{\gi}[1]{#1\mbox{{\rm -Ginj}}}

\newcommand{\opp}{^{\rm op}}
\newcommand{\otimesL}{\otimes^{\rm\mathbb L}}
\newcommand{\rHom}{{\rm\mathbb R}{\rm Hom}\,}
\newcommand{\pdim}{\pd}
\newcommand{\Hom}{{\rm Hom}}
\newcommand{\Coker}{{\rm Coker}}
\newcommand{ \Ker  }{{\rm Ker}}
\newcommand{ \Cone }{{\rm Con}}
\newcommand{ \Img  }{{\rm Im}}
\newcommand{\Ext}{{\rm Ext}}
\newcommand{\StHom}{{\rm \underline{Hom}}}
\newcommand{\StEnd}{{\rm \underline{End}}}

\newcommand{\KK}{I\!\!K}

\newcommand{\gm}{{\rm _{\Gamma_M}}}
\newcommand{\gmr}{{\rm _{\Gamma_M^R}}}

\def\vez{\varepsilon}\def\bz{\bigoplus}  \def\sz {\oplus}
\def\epa{\xrightarrow} \def\inja{\hookrightarrow}

\newcommand{\lra}{\longrightarrow}
\newcommand{\llra}{\longleftarrow}
\newcommand{\lraf}[1]{\stackrel{#1}{\lra}}
\newcommand{\llaf}[1]{\stackrel{#1}{\llra}}
\newcommand{\ra}{\rightarrow}
\newcommand{\dk}{{\rm dim_{_{k}}}}

\newcommand{\colim}{{\rm colim\, }}
\newcommand{\limt}{{\rm lim\, }}
\newcommand{\Add}{{\rm Add }}
\newcommand{\Prod}{{\rm Prod }}
\newcommand{\Tor}{{\rm Tor}}
\newcommand{\Cogen}{{\rm Cogen}}
\newcommand{\Tria}{{\rm Tria}}
\newcommand{\Loc}{{\rm Loc}}
\newcommand{\Coloc}{{\rm Coloc}}
\newcommand{\tria}{{\rm tria}}
\newcommand{\Con}{{\rm Con}}
\newcommand{\Thick}{{\rm Thick}}

{\Large \bf
\begin{center}
Virtually Gorenstein algebras of infinite dominant dimension
\end{center}}

\medskip
\medskip\centerline{\textbf{Hongxing Chen} and \textbf{Changchang Xi}$^*$}
\renewcommand{\thefootnote}{\alph{footnote}}
\setcounter{footnote}{-1} \footnote{ $^*$ Corresponding author.
Email: xicc@cnu.edu.cn; Fax: 0086 10 68903637.}
\renewcommand{\thefootnote}{\alph{footnote}}
\setcounter{footnote}{-1} \footnote{2020 Mathematics Subject
Classification: Primary 18G65, 16G10, 18G20; Secondary 16E65,16E35.}
\renewcommand{\thefootnote}{\alph{footnote}}
\setcounter{footnote}{-1} \footnote{Keywords: Dominant dimension; Nakayama conjecture; Self-orthogonal module; Virtually Gorenstein algebra.}

\begin{abstract}
Motivated by understanding the Nakayama conjecture which states that algebras of infinite dominant dimension should be self-injective, we study self-orthogonal modules with virtually Gorenstein endomorphism algebras and prove the following result: Given a finitely generated, self-orthogonal module over an Artin algebra with an orthogonal condition on its Nakayama translation, if its endomorphism algebra is virtually Gorenstein, then the module is projective. As a consequence, we re-obtain a recent result: the Nakayama conjecture holds true for the class of strongly Morita, virtually Gorenstein algebras. Finally, we show that virtually Gorenstein algebras can be constructed from Frobenius extensions.
\end{abstract}

\section{Introduction}\label{Introduction}
The dominant dimensions of algebras were introduced by Nakayama in 1958 (see \cite{Nakayama}) and have played an important role in the representation theory and homological algebra of finite-dimensional algebras. They have been studied intensively by Tachikawa, Morita, M\"{u}ller and many others (for example, see \cite{CFKKY,xc6, FK11,FK14,OI, mar, Muller,Nakayama,nrtz,Tac}).

\begin{Def}\label{dom-dim}
Let $A$ be an Artin algebra. The {\em dominant dimension} of $A$, denoted by  $\dm(A)$, is the largest natural number $n$
or $\infty$, such that, in a minimal injective coresolution
$$0\lra {_A}A\lra I^0\lra I^1\lra\cdots\lra I^{i-1}\lra I^i\lra\cdots$$
of the regular $A$-module $_AA$, all $I^i$ are projective for $0\leq i <n$.
\end{Def}

Doninant dimensions are closely related to self-orthogonal generator-cogenerators. Recall that a finitely generated module $M$ over an Artin algebra $A$ is said to be \emph{self-orthogonal} if $\Ext_A^i(M, M)=0$ for all $i\geq 1$; and is called a \emph{generator-cogenerator} if all indecomposable projective $A$-modules and indecomposable injective $A$-modules are isomorphic to direct summands of $M$. According to the Morita-Tachikawa correspondence \cite{Mo58a,Tac}, Artin algebras of dominant dimension at least $2$ are exactly the endomorphism algebras of generator-cogenerators. Moreover, M\"{u}ller showed that the endomorphism algebra $B$ of a generator-cogenerator $M$ has dominant dimension at least $n>1$ if and only if $\Ext_A^i(M,M)=0$ for all $1\leq i<n-1$ (see \cite[Lemma 3]{Muller}).

The extreme case $n=\infty$ involves the \emph{Nakayama conjecture} (see \cite{Nakayama}), one of the core problems in representation theory and homological algebra of finite-dimensional algebras (see \cite[p.409-410 ]{ars}):

\smallskip
\textbf{(NC)} If an Artin algebra has infinite dominant dimension, then it is self-injective.

\smallskip
This conjecture can be interpreted equivalently by self-orthogonal modules as follows \cite{Muller}:

\smallskip
\textbf{(NC-M)} A generator-cogenerator over an Artin algebra is projective  whenever it is self-orthogonal.

\smallskip
To understand the Nakayama conjecture, Tachikawa considered special self-orthogonal modules and divided the Nakayama conjecture into two conjectures, called \emph{Tachikawa's first and second conjectures} nowadays (see \cite[p.\! 115-116]{Tac}.

\smallskip
{\bf (TC1)} If an Artin algebra $A$  satisfies $\Ext^n_A(D(A),A)=0$ for all $n\geq 1$, then $A$ is self-injective, where $D$ is the usual duality of Artin algebra.

{\bf (TC2)}  Let $A$ be a self-injective Artin algebra and $M$ a finitely generated $A$-module. If $M$ is self-orthogonal, then $M$ is projective.

\smallskip
For a collection of all related conjectures and open problems, we refer to \cite[Conjcetures , p.409; open problems, p.411]{ars}. It is known in \cite{Tac} that (NC) holds if and only if both (TC1) and (TC2) hold.

Despite of efforts made in the past decades, all these conjectures still remain open in general.
Recently, some new advances on Tachikawa's second conjecture and the Nakayama conjecture have been made in \cite{xcf,xc3}. It is proved that Tachikawa's second conjecture for symmetric algebras is equivalent to saying that indecomposable symmetric algebras do not have any non-trivial stratifying ideals (see \cite[Theorem 1.1]{xcf}). Moreover, it is shown that the Nakayama conjecture holds for Gorenstein-Morita algebras introduced in \cite{xc3}. One of the main tools there to prove these results is recollements of certain ``nice'' triangulated categories such as Gorenstein stable categories or derived module categories of algebras.

In the present paper, we consider a self-orthogonal generator-cogenerator $M$ over an arbitrary Artin algebra, such that its Nakayama translation is orthogonal to $M$. If the endomorphism algebra of $M$ is virtually Gorenstein in the sense of Beligiannis and Reiten (see \cite[Chapter X, Definition 3.3]{BI}), then $M$ is projective. Our proof is based on an amazing characterization of virtually Gorenstein algebras in terms of contravariantly finite subcategories of module categories given in \cite{BK}. As a corollary of our main results, we re-obtain a recent result in \cite{xc3}: Strongly Morita, virtually Gorenstein algebras satisfy the Nakayama conjecture.

\medskip
To state our result more precisely, we introduce a few notions and notation.

Unless stated otherwise, all algebras considered are Artin algebras over a fixed commutative Artin ring, and all
modules are left modules, unless stated otherwise.

Let $A$ be an algebra. We denote by $A$-Mod (or $A$-mod) the category of all (or finitely generated) $A$-modules, and by $\nu_A$ the Nakayama functor $D(A)\otimes_A-$, where $D$ stands for the usual duality over Artin algebras universally.

Following \cite[Definition 8.1]{Bel}, an algebra $A$ is said to be \emph{virtually Gorenstein} provided that for each $A$-module $X\in A\Modcat$, the functor $\Ext_A^i(X,-)$ vanishes for all $i\geq 1$ on all Gorenstein injective $A$-modules in $A\Modcat$ if and only if the functor $\Ext_A^i(-, X)$ vanishes for all $i\geq 1$ on all Gorenstein projective $A$-modules in $A\Modcat$. The class of virtually Gorenstein algebras contains Gorenstein algebras and algebras of finite representation type, and is closed under taking derived equivalences and stably equivalences of Morita type (see \cite{Bel,Bel2, BK}). Moreover, it was shown in \cite{Bel} that virtually Gorenstein algebras satisfy the Gorenstein symmetric conjecture (see \cite[Conjecture (13), p.410]{ars} for the statement). Note, however, that not all algebras are virtually Gorenstein (see \cite{BK} for a counterexample). As a generalization of virtually Gorenstein algebras, the class of compactly Gorenstein algebras is introduced in \cite[Section 1.2]{xc3}. We conjecture that all Artin algebras should be compactly Gorenstein.

One of our main results is a combination of the Nakayama conjecture and Tachikawa's first conjecture on virtually Gorenstein algebras (see Remark \ref{TC1} for the first conjecture).

\begin{Theo}\label{DM+EXT}
Suppose that $A$ is an algebra with $\dm(A)=\infty$ and $\Ext_A^n(D(A), A)=0$ for all $n\geq 1$.
If $A$ is virtually Gorenstein, then $A$ is self-injective.
\end{Theo}

Theorem \ref{DM+EXT} will be used to prove our next result about self-orthogonal modules.

\begin{Theo}\label{Main result}
Let $A$ be an algebra, and let $M$ be a finitely generated, generator-cogenerator for $A\modcat$. Suppose $\Ext_A^n(M\oplus \nu_A(M), M)=0$ for all $n\geq 1$. If the endomorphism algebra of the $A$-module $M$ is virtually Gorenstein, then $M$ is a projective $A$-module.
\end{Theo}

The above results reveal a close relation between the Nakayama conjecture and virtually Gorenstein algebras. A direct consequence of Theorem \ref{Main result} is the following corollary which includes the case that $A$ is a symmetric algebra. In this case, the Nakayama functor is the identity functor.

\begin{Koro}\label{Nakayama stable}
Let $A$ be an algebra and let $M$ be a finitely generated, self-orthogonal $A$-module which is a generator-cogenerator with $\nu_A(M)\in\add(_AM)$. If the endomorphism algebra of $_AM$ is virtually Gorenstein, then $M$ is a projective $A$-module.
\end{Koro}

Recall from \cite[Section 1.2]{xc3} that \emph{strongly Morita} algebras are, by definition, the endomorphism algebras of those generators $M$ over a self-injective algebra $A$ such that $\add(_AM)=\add(\nu_A(M))$. The class of strongly Morita algebras contains \emph{gendo-symmetric algebras} that are the endomorphism algebras of generators over symmetric algebras (see \cite{FK11,FK14}).

Now, we apply Corollary \ref{Nakayama stable} to strongly Morita algebras and give a completely \emph{different} proof of the following result which is a special case of \cite[Corollary 1.4]{xc3}.

\begin{Koro}\label{Gendo-symmetric}
Let $A$ be a strongly Morita, virtually Gorenstein Artin algebra. If $A$ has infinite dominant dimension, then it is self-injective.
\end{Koro}

{\bf Proof.} Let $A$ be a strongly Morita algebra. Then $A=\End_{\Lambda}(M)$, where $\Lambda$ is a self-injective algebra and $M\in \Lambda\modcat$ is a generator with $\add(_{\Lambda}M)=\add(\nu_{\Lambda}(M))$. Suppose $\dm(A)=\infty$. Then $\Ext_{\Lambda}^n(M, M)=0$ for all $n\geq 1$ by \cite[Lemma 3]{Muller}. Since $A$ is virtually Gorenstein, the $\Lambda$-module $M$ is projective by Corollary \ref{Nakayama stable}. This implies that $A$ is Morita equivalent to $\Lambda$, and thus self-injective because Morita equivalences preserve self-injective algebras.  $\square$

\smallskip
Finally, we point out that virtually Gorenstein algebras can be obtained from Frobenius extensions. For details, we refer the reader to Proposition \ref{frob-ext}.

\bigskip
The contents of this paper are sketched as follows. In Section $2$, we fix some notation and recall the definitions of contravariantly or covariantly finite subcategories as well as two relevant theorems. In Section $3$, we first give some properties of algebras of infinite dominant dimensions
(Lemma \ref{Thick}) and then show Theorem \ref{DM+EXT}. Subsequently, we apply Theorem \ref{DM+EXT} to show Theorem \ref{Main result}. Finally, we show that Frobenius extensions provide a way to get new virtually Gorenstein algebras from given ones.

\section{Preliminaries\label{sect2}}
In this section, we briefly recall some definitions and notation used in this paper.

Let $\mathcal C$ be an additive category.

Let $X$ be an object in $\mathcal{C}$. We denote by $\add(X)$ the full subcategory of $\mathcal{C}$ consisting of all direct summands of
finite coproducts of copies of $M$. If $\mathcal{C}$ admits small coproducts (that is, coproducts indexed over sets exist in
${\mathcal C}$), then we denote by $\Add(X)$ the full subcategory of $\mathcal{C}$ consisting of all direct summands of small coproducts
of copies of $X$. Dually, if $\mathcal{C}$ admits products, then $\Prod(X)$ denotes the full subcategory of $\mathcal{C}$
consisting of all direct summands of products of copies of $X$.

Let $\mathcal{B}$ be a full subcategory of $\mathcal{C}$. A morphism $f: X\to Y$ in $\mathcal{C}$ is called a
\emph{right $\mathcal{B}$-approximation} of $Y$ provided that $X\in \mathcal{B}$ and
$\Hom_\mathcal{C}(B, f): \Hom_\mathcal{C}(B, X)\to \Hom_\mathcal{C}(B, Y)$ is surjective for any $B\in\mathcal{B}$.
If each object of $\mathcal{C}$ admits a right $\mathcal{B}$-approximation, then $\mathcal{B}$ is said to be \emph{contravariantly finite}.
Dually, we can define left approximations of objects and \emph{covariantly finite} subcategories in $\mathcal{C}$.

Let $\mathcal{C}$ be an abelian category. The category $\mathcal{B}$ is called a \emph {thick subcategory} of $\mathcal{C}$
if it is closed under direct summands in $\mathcal{C}$ and has
the \emph {two out of three property}: for any short exact sequence
$0\to X\to Y\to Z\to 0$ in $\mathcal{C}$ with two terms in $\mathcal{B}$, the third term belongs to $\mathcal{B}$ as well.
For a class $\mathcal{S}$ of objects in $\mathcal{C}$, we denote by $\Thick(\mathcal{S})$ the smallest thick subcategory of $\mathcal{C}$ containing $\mathcal{S}$. When $\mathcal{C}$ has enough projective objects,
$\mathcal{B}$ is called a \emph{resolving subcategory} of $\mathcal{C}$ if $\mathcal{B}$ contains all projective objects of $\mathcal{C}$ and is closed under extensions and kernels of epimorphisms in $\mathcal{C}$ (see \cite[Section 3]{AR2}).

Let $A$ be an Artin algebra. Recall that $A\Modcat$ (respectively, $A\modcat$) denotes the category of all (respectively, finitely generated) left $A$-modules.
Let $\Omega_A$ and $\Omega_A^{-}$ stand for the usual syzygy and cosyzygy functors over $A\Modcat$, respectively.
For a class $\mathcal{S}$ of objects in $A\Modcat$, we denote by $\mathcal{S}^{\bot}$ (respectively, ${^{\bot}}\mathcal{S}$) the full subcategory of $A\Modcat$ consisting of modules $X$ such that $\Ext_A^n(S, X)=0$ (respectively, $\Ext_A^n(X,S)=0$) for all $S\in\mathcal{S}$ and $n\geq 1$.

Let $\mathcal{U}$ and $\mathcal{V}$ be full subcategories of $A\Modcat$ closed under isomorphisms.
Denote by $\mathcal{U}\widehat{\oplus}\mathcal{V}$ the full subcategory of $\Modcat$ which consists of all modules $X$ such that
$X\oplus W\simeq U\oplus V$, where $W\in\mathcal{U}\cap\mathcal{V}$, $U\in\mathcal{U}$ and $V\in\mathcal{V}$.
Note that if $\mathcal{U}\subseteq A\modcat$ and $\mathcal{V}\subseteq A\modcat$ are closed under direct summands, then $X\in \mathcal{U}\widehat{\oplus}\mathcal{V}$ if and only if
$X\simeq U\oplus V$ with $U\in\mathcal{U}$ and $V\in\mathcal{V}$. In this case, we simply write $\mathcal{U}\oplus\mathcal{V}$ for $\mathcal{U}\widehat{\oplus}\mathcal{V}$;
in other words, $\mathcal{U}\oplus\mathcal{V}:=\{X\in A\modcat\mid X\simeq U\oplus V,\; U\in\mathcal{U}, V\in\mathcal{V}\}$.

Let $\Prj{A}$ and $\Inj{A}$ (respectively, $\prj{A}$ and $\inj{A}$) be the full subcategories of $A\Modcat$ consisting of (respectively, finitely generated)
projective and injective $A$-modules, respectively.
As usual, the projective and injective dimensions of an $A$-module $X$ are denoted by $\pd(X)$ and $\id(X)$, respectively.
Let $$\mathscr{P}^{<\infty}(A):=\{X\in B\Modcat\mid \pd(X)<\infty\}\quad \mbox{and}\quad \mathscr{I}^{<\infty}(A):=\{X\in B\Modcat\mid \id(X)<\infty\}.$$
They are thick subcategories of $A\Modcat$. Their restrictions to finitely generated modules are denoted by
$$\mathscr{P}^{<\infty}_{{\rm fg}}(A):=\mathscr{P}^{<\infty}(A)\cap A\modcat\quad \mbox{and}\quad \mathscr{I}^{<\infty }_{{\rm fg}}(A):=\mathscr{I}^{<\infty}(A)\cap A\modcat.$$ Then $\mathscr{P}^{<\infty}(A)$ is a resolving subcategory, and $\mathscr{P}^{<\infty}_{{\rm fg}}(A)\cup\mathscr{I}^{<\infty }_{{\rm fg}}(A)\subseteq \Thick(\prj{A}\cup \inj{A})\subseteq A\modcat$.

\smallskip
As a preparation for showing Theorem \ref{Main result}, we need the following two important results. The first one characterizes virtually Gorenstein algebras.

\begin{Theo} {\rm \cite[Theorem 1]{BK}}\label{BK}
The following are equivalent for an Artin algebra $A$.

$(1)$ The algebra $A$ is virtually Gorenstein.

$(2)$ The subcategory $\Thick(\prj{A}\cup \inj{A})$ of $A\modcat$ is contravariantly finite.

$(3)$ The subcategory $\Thick(\prj{A}\cup \inj{A})$ of $A\modcat$ is covariantly finite.
\end{Theo}

The next result describes modules in a resolving, contravariantly finite subcategory.

\begin{Theo} {\rm \cite[Proposition 3.8]{AR2}}\label{AUR}
Let $A$ be an algebra.
Suppose $\mathscr{X}$ is a resolving, contravariantly finite subcategory
of $A\modcat$. Let $S_1, S_2,\cdots, S_t$ be a complete set of nonisomorphic simple
$A$-modules and let $f_i: X_i\to S_i$ be a minimal right $\mathscr{X}$-approximation of $S_i$
for $1\leq i\leq t$. Then the modules in $\mathscr{X}$ consist of the summands of modules $M$
with the property that there is a finite filtration $M=M_0\supset M_1\supset \cdots \supset M_n=0$ such that, for each $0\leq i\leq n-1$, we have $M_i/M_{i+1}\simeq X_j$ for some $j\in \{1, 2, \cdots, t\}$.
\end{Theo}

\section{Algebras of infinite dominant dimension}

In this section, we are concentrated on algebras of infinite dominant dimension. These algebras have the following property.

\begin{Lem}\label{Thick}
Let $B$ be an algebra of infinite dominant dimension. Then the following hold true.

$(1)$ There is a finitely generated $B$-module $E$ such that $\Add(E)=\Prj{B} \cap \Inj{B}=\mathscr{P}^{<\infty}(B)\cap\mathscr{I}^{<\infty}(B).$

$(2)$ Suppose $\Ext_B^n(D(B), B)=0$ for all $n\geq 1$. Then $\Thick(\Prj{B}\cup \Inj{B})=\mathscr{P}^{<\infty}(B)\widehat{\oplus}\mathscr{I}^{<\infty}(B)$ and $\Thick(\prj{B}\cup \inj{B})=\mathscr{P}^{<\infty}_{{\rm fg}}(B)\oplus\mathscr{I}^{<\infty}_{{\rm fg}}(B)$.
\end{Lem}
{\it Proof.} Since $B$ is an Artin algebra, it is known that each finitely generated $B$-module $M$ satisfies $\Add(M)=\Prod(M)$ (for example, see \cite[Lemma 1.2]{KS}). This property will be used freely in our proof.

$(1)$ Let $\mathscr{E}:=\Prj{B}\cap\Inj{B}$ be the category of projective-injective $B$-modules, and let $E\in B\modcat$ such that $\add(E)=\prj{B}\cap\inj{B}$.
Then $\Add(E)\subseteq\mathscr{E}$. Since $\dm(B)\ge 1$, the injective envelope of ${_B}B$ belongs to $\add(E)$. It follows that each projective $B$-module can be embedded into a module in $\Add(E)$, and therefore $\mathscr{E}\subseteq \Add(E)$ by the splitting property of injective modules.
Thus $\mathscr{E}=\Add(E)$.

Let $\mathscr{E}\mbox{-dim}_{\infty}(B)$ (respectively, $\mathscr{E}\mbox{-dim}^{\infty}(B)$) be the full subcategory of $B\Modcat$ consisting of all modules $X$ such that there is a long exact sequence of $B$-modules
$$\cdots\lra X_2\lra X_1\lra X_0\lra X\lra 0 \quad (\mbox{respectively,}\; 0\lra X\lra X_0\lra X_1\lra X_2\lra \cdots )$$ with $X_i\in\mathscr{E}$ for all $i\geq 0$. As $\mathscr{E}$ consists of all projective-injective $B$-modules and is a thick subcategory of $B\Modcat$, we can show that both $\mathscr{E}\mbox{-dim}_{\infty}(B)$ and $\mathscr{E}\mbox{-dim}^{\infty}(B)$ are thick subcategories of $B\Modcat$ containing $\mathscr{E}$. Moreover, since $\mathscr{E}=\Add(E)=\Prod(E)$, the categories $\mathscr{E}\mbox{-dim}_{\infty}(B)$ and $\mathscr{E}\mbox{-dim}^{\infty}(B)$ are closed under direct sums and products in $B\Modcat$. Since $\dm(B)=\infty$, we have ${_B}B\in \mathscr{E}\mbox{-dim}^{\infty}(B)$. Note that $\dm(B\opp)=\dm(B)$ by \cite[Theorem 4]{Muller}. Thus $D(B_B)\in \mathscr{E}\mbox{-dim}_{\infty}(B)$. Consequently,
$$ (\dag)\qquad
\mathscr{P}^{<\infty}(B)\subseteq \mathscr{E}\mbox{-dim}^{\infty}(B)\quad\mbox{and}\quad \mathscr{I}^{<\infty}(B)\subseteq\mathscr{E}\mbox{-dim}_{\infty}(B).
$$ Since $\mathscr{E}\mbox{-dim}^{\infty}(B)\cap\mathscr{I}^{<\infty}(B)=\mathscr{E}$, we obtain
$\mathscr{P}^{<\infty}(B)\cap\mathscr{I}^{<\infty}(B)=\mathscr{E}$. This shows $(1)$.

$(2)$ Let $\mathscr{C}:=\mathscr{P}^{<\infty}(B)\widehat{\oplus}\mathscr{I}^{<\infty}(B)$. Since $\mathscr{P}^{<\infty}(B)=\Thick(\Prj{B})$ and $\mathscr{I}^{<\infty}(B)=\Thick(\Inj{B})$, we have $\mathscr{C}\subseteq\Thick(\Prj{B}\cup \Inj{B})$. To show the converse inclusion,
it suffices to show that $\mathscr{C}$ is a thick subcategory of $B\Modcat$. However, this will be done by the following three steps.

{\bf Step 1.} We show that $\mathscr{C}$ is closed under extensions in $B\Modcat$.

In fact, by the assumption of $(2)$, $D(B_B)\in {^{\bot}}B$. It follows from $\Inj{B}=\Prod(D(B_B))=\Add(D(B_B))$ that $\Inj{B}\subseteq{^{\bot}}B$.
Note that the category ${^{\bot}}B$ is always closed under kernels of surjections in $B\Modcat$. This implies $\mathscr{I}^{<\infty}(B)\subseteq {^{\bot}}B$, or equivalently,
${_B}B\in\mathscr{I}^{<\infty}(B)^{\bot}$. Since $\Prj{B}=\Prod({_B}B)$ and $\mathscr{I}^{<\infty}(B)^{\bot}$ is closed under cokernels of injections in $B\Modcat$, we have $\mathscr{P}^{<\infty}(B)\subseteq \mathscr{I}^{<\infty}(B)^{\bot}$. Further, we show $\mathscr{I}^{<\infty}(B)\subseteq\mathscr{P}^{<\infty}(B)^{\bot}$. In fact, for any $X\in\mathscr{P}^{<\infty}(B)$ and for any $n\in\mathbb{N}$, it follows from
$\mathscr{P}^{<\infty}(B)\subseteq\mathscr{E}\mbox{-dim}^{\infty}(B)$ (see the inclusion in (\dag)) that there are $B$-modules $X_n\in B\Modcat$ and $Q_n\in\Add(_BB)$ such that
$X\simeq \Omega_B^n(X_n)\oplus Q_n$. For $Y\in B\Modcat$ and $m\geq 1$, it is clear that
$$\Ext_B^m(X,Y)\simeq\Ext_B^m(\Omega_B^n(X_n)\oplus Q_n,Y)\simeq \Ext_B^{m+n}(X_n,Y)\simeq\Ext_B^m(X_n,\Omega_B^{-n}(Y)).$$
Thus $\Ext_B^m(X,Y)=0$ for $m\geq 1$ if $Y\in\mathscr{I}^{<\infty}(B)$. This shows $\mathscr{I}^{<\infty}(B)\subseteq\mathscr{P}^{<\infty}(B)^{\bot}$.

To complete the proof of Step 1, we apply $\mathscr{P}^{<\infty}(B)\subseteq \mathscr{I}^{<\infty}(B)^{\bot}$ and $\mathscr{I}^{<\infty}(B)\subseteq\mathscr{P}^{<\infty}(B)^{\bot}$ to show the following fact $(\ast)$.

\smallskip
$(\ast)\;$ Each exact sequence $0\to U_1\oplus V_1\to W\to U_2\oplus V_2\to 0$ of $B$-modules with $U_i\in\mathscr{P}^{<\infty}(B)$ and $V_i\in\mathscr{I}^{<\infty}(B)$ for $i=1,2$, is isomorphic to a direct sum of two exact sequences $0\to U_1\to W_1\to U_2\to 0$ and $0\to V_1\to W_2\to V_2\to 0$ in $B\Modcat$. In particular, $W\simeq W_1\oplus W_2$ with $W_1\in\mathscr{P}^{<\infty}(B)$ and $W_2\in \mathscr{I}^{<\infty}(B)$.

\smallskip
Indeed, it follows from $\mathscr{I}^{<\infty}(B)\subseteq\mathscr{P}^{<\infty}(B)^{\bot}$ that $\Ext_B^1(U_2,  V_1)=0$, and from $\mathscr{P}^{<\infty}(B)\subseteq \mathscr{I}^{<\infty}(B)^{\bot}$ that $\Ext_B^1(V_2,  U_1)=0$. Hence, by the finite additivity of the bifunctor $\Ext^i_B(-,-)$, we have the following isomorphism of abelian groups
$$(\ddag)\qquad \Ext_B^1(U_2\oplus V_2, U_1\oplus V_1)\simeq \Ext_B^1(U_2, U_1)\oplus\Ext_B^1(V_2, V_1).$$
Note that, for a pair $(U,V)$ of $B$-modules, $\Ext_B^1(V,U)$ can be interpreted as the abelian group of the equivalence classes of short exact sequences $0\to U\to E\to V\to 0$ in $B\Modcat$ (for example, see \cite[I. Theorem 5.4]{ars}).
Thus $(\ast)$ follows from interpreting $(\ddag)$ as short exact sequences.

Now, it follows from $(\ast)$ and $\Add(E)=\mathscr{P}^{<\infty}(B)\cap\mathscr{I}^{<\infty}(B)$ that
$\mathscr{C}$ is closed under extensions in $B\Modcat$.

{\bf Step 2.} We show that $\mathscr{C}$ is closed under direct summands in $B\Modcat$. Alternately, we show that if $X\oplus Y\simeq U\oplus V$ in $B\Modcat$ with $U\in \mathscr{P}^{<\infty}(B)$ and $V\in\mathscr{I}^{<\infty}(B)$, then $X\in\mathscr{C}$.

Let $n:=\pd(U)+1<\infty$. Then $\Omega_B^n(U)=0$ which leads to  $\Omega_B^n(X)\oplus\Omega_B^n(Y)\simeq\Omega_B^n(V)$.
According to $V\in\mathscr{I}^{<\infty}(B)\subseteq\mathscr{E}\mbox{-dim}_{\infty}(B)$, there is an  exact sequence
$$0\lra \Omega_B^n(V)\lra P_{n-1}\lra \cdots\lra P_1\lra P_0\lra V\lra 0,$$ where $P_i$ lies in $\Add(E)$ for  $0\leq i\leq n-1$.
Let
$$0\lra \Omega_B^n(X)\lra I^0\lra I^1\lra \cdots \lra I^{n-1}\lra I^n\lra\cdots $$ be a minimal injective coresolution of $\Omega_B^n(X)$. Since all $P_i$ are injective and $\Omega_B^n(X)$ is a direct summand of $\Omega_B^n(V)$, the injective module $I^j$ is a direct summand of $P_{n-1-j}$ for $0\leq j\leq n-1$, and $\Omega_B^{-n}(\Omega_B^n(X))$ is a direct summand of $V$. In particular, $I^j$ is projective and $\Omega_B^{-n}(\Omega_B^n(X))\in\mathscr{I}^{<\infty}(B)$. Now, we can construct the following exact commutative diagram
$$
\xymatrix{
0\ar[r] &\Omega^n_B(X)\ar[r]\ar@{=}[d] & Q_{n-1}\ar[r]\ar@{-->}[d]_-{f_0}&\cdots \ar[r]& Q_1\ar[r]\ar@{-->}[d]_-{f_{n-2}}& Q_0\ar[r] \ar@{-->}[d]_-{f_{n-1}} &X\ar[r]\ar@{-->}[d]_-{\varepsilon_X} & 0\\
0\ar[r] &\Omega^n_B(X)\ar[r] & I^0\ar[r] &\cdots \ar[r]& I^{n-2}\ar[r]& I^{n-1}\ar[r] & \Omega_B^{-n}(\Omega_B^n(X))\ar[r] & 0
}$$
in which the first row arises from a minimal projective resolution of $X$ and vertical maps are induced from the identity map of $\Omega^n_B(X)$.
Taking the mapping cone of the quasi-isomorphism $(f_0, f_1,\cdots, f_{n-1}, \varepsilon_X)$ yields a long exact sequence
$$
0\lra Q_{n-1}\lra  Q_{n-2}\oplus I^0\lra \cdots\lra Q_0\oplus I^{n-2}\lraf{h} X\oplus I^{n-1}\lra \Omega_B^{-n}(\Omega_B^n(X))\lra 0.
$$
This implies $L:=\Img(h)\in\mathscr{P}^{<\infty}(B)$ because $Q_i$ and $I^i$ are projective modules for all $0\leq i\leq n-1$. Since $\Omega_B^{-n}(\Omega_B^n(X))\in\mathscr{I}^{<\infty}(B)$ and $\mathscr{P}^{<\infty}(B)\subseteq \mathscr{I}^{<\infty}(B)^{\bot}$, we have
$\Ext_B^1\big(\Omega_B^{-n}(\Omega_B^n(X)), L\big)=0$. It follows that
$$X\oplus I^{n-1}\simeq L\oplus \Omega_B^{-n}(\Omega_B^n(X)).$$ Clearly, $I^{n-1}\in\Add(E)=\mathscr{P}^{<\infty}(B)\cap\mathscr{I}^{<\infty}(B)$, $L\in\mathscr{P}^{<\infty}(B)$ and $\Omega_B^{-n}(\Omega_B^n(X))\in\mathscr{I}^{<\infty}(B)$. Thus $X\in\mathscr{C}$.

Remark that if $X$ is finitely generated, then all the modules in the above commutative diagram are finitely generated. In this situation, $I^{n-1}\in\add(E)$, $L\in \mathscr{P}^{<\infty}_{{\rm fg}}(B)$ and $\Omega_B^{-n}(\Omega_B^n(X))\in\mathscr{I}^{<\infty}_{{\rm fg}}(B)$.

{\bf Step 3.} We show that $\mathscr{C}$ is closed under kernels of surjective homomorphisms, and cokernels of injective homomorphisms in $B\Modcat$.

Actually, since $\mathscr{E}$ consists of projective-injective modules and $\mathscr{I}^{<\infty}(B)\subseteq\mathscr{E}\mbox{-dim}_{\infty}(B)$, there holds $\Omega_B(\mathscr{I}^{<\infty}(B))\subseteq \mathscr{I}^{<\infty}(B)$. Clearly, $\Omega_B(\mathscr{P}^{<\infty}(B))\subseteq \mathscr{P}^{<\infty}(B)$. Thus $\Omega_B(\mathscr{C})\subseteq \mathscr{C}$. Dually, $\Omega_B^{-}(\mathscr{C})\subseteq\mathscr{C}$.

Let $(\delta): 0\to X\to Y\to Z\to 0$ be an exact sequence in $B\Modcat$. Then there are
two relevant exact sequences in $B\Modcat$:
$$
(\delta_1):\;\; 0\lra\Omega_B(Z)\lra X\oplus P_Z\lra Y\lra 0\quad and \quad (\delta_2):\;\; 0\lra Y\lra Z\oplus I_X\lra \Omega_B^-(X)\lra 0,
$$ where $P_Z$ is a projective cover of $Z$ and $I_X$ is an injective envelop of $X$.
To show Step $3$, we consider the following two cases:

$(a)$ Suppose both $Y$ and $Z$ lie in $\mathscr{C}$. Then $\Omega_B(Z)\in\mathscr{C}$. It follows from $(\delta_1)$ and Step $1$ that $X\oplus P_Z\in\mathscr{C}$. Thus $X\in\mathscr{C}$ by Step $2$.
This means that $\mathscr{C}$ is closed under kernels of surjections in $B\Modcat$.

$(b)$ Suppose both $X$ and $Y$ lie in $\mathscr{C}$. Then $\Omega_B^-(X)\in\mathscr{C}$. It follows from $(\delta_2)$ and Step $1$ that $Z\oplus I_X\in\mathscr{C}$. Thus $Z\in\mathscr{C}$ by Step $2$.
Hence $\mathscr{C}$ is closed under cokernels of injections in $B\Modcat$. This completes Step 3.

Thus $\mathscr{C}$ is a thick subcategory of $B\Modcat$, and $\mathscr{C}=\Thick(\Prj{B}\cup \Inj{B})$.
Similarly, by considering finitely generated $B$-modules, we can prove the equality
$$\Thick(\prj{B}\cup \inj{B})=\mathscr{P}^{<\infty}_{{\rm fg}}(B)\widehat{\oplus}\mathscr{I}^{<\infty}_{{\rm fg}}(B).$$ Since both $\mathscr{P}^{<\infty}_{{\rm fg}}(B)$ and $\mathscr{I}^{<\infty}_{{\rm fg}}(B)$ are closed under direct summands in $B\modcat$, we have $\mathscr{P}^{<\infty}_{{\rm fg}}(B)\widehat{\oplus}\mathscr{I}^{<\infty}_{{\rm fg}}(B)$ = $\mathscr{C}\cap B\modcat=\mathscr{P}^{<\infty}_{{\rm fg}}(B)\oplus\mathscr{I}^{<\infty}_{{\rm fg}}(B)$. This shows $(2)$.
$\square$

\begin{Rem}\label{TC1}
Given a finite-dimensional $k$-algebra $B$ over a field $k$, the condition $\Ext_B^n(D(B), B)=0$ for all integers $n\geq 1$ in Lemma \ref{Thick} is equivalent to saying that the \emph{minimal} self-orthogonal generator-cogenerator $B\oplus D(B)$ for $B\modcat$ is self-orthogonal. This is also related to \emph{Tachikawa's first conjecture}:
If $\Ext_{B\otimes_kB\opp}^n(B, B\otimes_kB)=0$ for all $n\geq 1$, then $B$ is self-injective (see \cite[p.\! 115]{Tac}).

In fact, there are isomorphisms of $k$-modules for all $n\geq 1$:
$$
\Ext_B^n(B\oplus D(B), B\oplus D(B))\simeq \Ext_B^n(D(B), B)\simeq \Ext_{B\otimes_kB\opp}^n(B, B\otimes_kB),
$$(see \cite[p.\! 114]{Tac} for the last isomorphism).
\end{Rem}

\medskip
{\bf Proof of Theorem \ref{DM+EXT}}.
Let $A$ be a virtually Gorenstein algebra, and let $\mathscr{X}=\Thick(\prj{A}\cup \inj{A})$. By Theorem \ref{BK}, $ \mathscr{X} $ is contravariantly finite in $A\modcat$. In other words, each finitely generated $A$-module $X$ has a minimal right $\mathscr{X}$-approximation $W_X\to X$, that is, $W_X\in\mathscr{X}$ and the induced map $\Hom_A(W, W_X)\to \Hom_A(W, X)$ is surjective for any $W\in\mathscr{X}$. Let $S_1, \cdots, S_m$ be a complete set of nonisomorphic simple $A$-modules, and let $W_i\to S_i$ be a minimal right $\mathscr{X}$-approximation of $S_i$ for $i=1,\cdots, m$. Since $\mathscr{X}=\mathscr{P}^{<\infty}_{{\rm fg}}(A)\oplus\mathscr{I}^{<\infty}_{{\rm fg}}(A)$ by Lemma \ref{Thick}(2), we have $W_i\simeq U_i\oplus V_i$ for some $U_i\in\mathscr{P}^{<\infty}_{{\rm fg}}(A)$ and $V_i\in\mathscr{I}^{<\infty}_{{\rm fg}}(A)$. Clearly, $\mathscr{X}$ is a thick subcategory of $A\modcat$ and contains all finitely generated projective $A$-modules. In particular, $\mathscr{X}$ is closed under extensions and kernels of surjections in $A\modcat$.
Thus $\mathscr{X}$ is a resolving subcategory of $A\modcat$. By Theorem \ref{AUR}, $\mathscr{X}$ consists of the direct summands of modules $X$ with a filtration of finite length $n$:
$$X=X_0\supset X_1\supset \cdots \supset X_n=0$$ such that, for each $j=0, \cdots, n-1$, there is an isomorphism $X_j/X_{j+1}\simeq W_{\ell_j}$ for some $\ell_j\in\{1,\cdots, m\}$.
Now, we fix such an $A$-module $X$. Since $W_{\ell_j}\simeq U_{\ell_j}\oplus V_{\ell_j}$ with $U_{\ell_j}\in\mathscr{P}^{<\infty}_{{\rm fg}}(A)$ and $V_{\ell_j}\in\mathscr{I}^{<\infty}_{{\rm fg}}(A)$, we see from  $(\ast)$ in the proof of Lemma \ref{Thick}(2) that there are finitely generated $A$-modules $Y$ and $Z$ with filtration of finite length:
$$Y=Y_0\supseteq Y_1\supseteq \cdots \supseteq Y_n=0\quad \mbox{and}\quad Z=Z_0\supseteq Z_1\supseteq \cdots \supseteq Z_n=0$$
such that for $0\leq j\leq n-1\in\mathbb{N}$,
$$X_j\simeq Y_j\oplus Z_j, \quad Y_j/Y_{j+1}\simeq U_{\ell_j}, \quad \mbox{and}\quad Z_j/Z_{j+1}\simeq V_{\ell_j}.$$ In particular, $X\simeq Y\oplus Z$ with $Y\in\mathscr{P}^{<\infty}_{{\rm fg}}(A)$ and $Z\in \mathscr{I}^{<\infty}_{{\rm fg}}(A)$. Set
$$s:=\max\{\pd(U_i)\mid 1\leq i\leq m\}\quad\mbox{and}\quad t:=\max\{\id(V_i)\mid 1\leq i\leq m\}.$$
Then $\pd(Y)\leq s$ and $\id(Z)\leq t$. Let $N$ be an indecomposable direct summand of $X$. Then $N$ is isomorphic to a direct summand of either $Y$ or $Z$. Thus $\pd(N)\leq s$ or $\id(N)\leq t$. Consequently, each indecomposable module in $\mathscr{X}$ has either projective dimension at most $s$ or injective dimension at most $t$. Now, let $T\in \mathscr{P}^{<\infty}_{{\rm fg}}(A)$. Then we can write
$T=\bigoplus_{0\leq i\leq u\in\mathbb{N}}T_i$ as a direct sum of indecomposable (finitely generated) $A$-modules $T_i$. Then either $\pd(T_i)\leq s$ or $\id(T_i)\leq t$. Recall that $\mathscr{P}^{<\infty}(A)\cap\mathscr{I}^{<\infty}(A)$ consists of projective-injective $A$-modules by Lemma \ref{Thick}(1). This implies that if $\id(T_i)\leq t$, then $\pd(T_i)=0$. Thus $\pd(T)=\max\{\pd(T_i)\mid 0\leq i\leq u\}\leq s$.

Since $\dm(A)=\infty$, the minimal injective coresolution $$0\lra {}_AA\lra I^0\lra I^1\lra\cdots \lra I^n\lra\cdots $$ of the module ${_A}A$ has all terms $I^i$ being projective-injective. This implies $\pd(\Omega_A^{-s-1}(A))\leq s+1<\infty$ and $\Omega_A^{-s-1}(A)\in \mathscr{P}^{<\infty}_{{\rm fg}}(A)$. Hence $\pd\big(\Omega_A^{-s}(\Omega_A^-(A))\big)=\pd(\Omega_A^{-s-1}(A))\leq s$. It follows that $\Omega_A^{-}(A)$ is projective, and therefore $I^0\simeq {}_AA\oplus \Omega_A^-(A)$. Thus ${_A}A$ is injective, as desired.  $\square$

\medskip
To show Theorem \ref{Main result}, we need the following result.

\begin{Lem}\label{WG}
Let $\Lambda$ be an algebra, $M$ a finitely generated $\Lambda$-module and $A$ the endomorphism algebra of $_{\Lambda}M$. Suppose that $_{\Lambda}M$ is a generator-cogenerator. Then $\dm(A)=\infty$ and $\Ext_A^n(D(A), A)=0$ for all $n\geq 1$ if and only if $\Ext_{\Lambda}^i(M\oplus \nu_{\Lambda}(M), M)=0$ for all $i\geq 1$.
\end{Lem}

{\it Proof.} By M\"{u}ller's theorem on dominant dimension (see \cite[Lemma 3]{Muller}), $\dm(A)=\infty$ if and only if $\Ext_{\Lambda}^i(M, M)=0$ for all $i\geq 1$. Now, we assume $\dm(A)=\infty$.

Let $$0\lra {}_{\Lambda}M\lra I_0\lra I_1\lra \cdots\lra I_n\lra \cdots$$ be a minimal injective coresolution of $_{\Lambda}M$. Note that $M$ is naturally a $\Lambda$-$A$-bimodule. Applying $\Hom_{\Lambda}(M,-)$ to this coresolution yields a long exact sequence of $A$-modules:
$$0\lra {_A}A\lra \Hom_{\Lambda}(M, I_0)\lra \Hom_{\Lambda}(M, I_1)\lra\cdots\lra\Hom_{\Lambda}(M, I_n)\lra\cdots$$
where $\Hom_{\Lambda}(M, I_n)$ are projective-injective for all $n\geq 0$.  In particular, this sequence is an injective coresolution of ${_A}A$. Now, we apply $\Hom_A(D(A),-)$ to the sequence and obtain a complex of $A$-modules:
$$
0\to\Hom_A(D(A), A)\to \Hom_A(D(A), \Hom_{\Lambda}(M, I_0))\to\cdots\to \Hom_A(D(A), \Hom_{\Lambda}(M, I_n))\to\cdots
$$
which is, by adjoint isomorphism, isomorphic to the complex
$$
(\ast\ast)\quad 0\to\Hom_{\Lambda}(M\otimes_AD(A), M)\to \Hom_{\Lambda}(M\otimes_AD(A), I_0)\to\cdots\to \Hom_{\Lambda}(M\otimes_AD(A), I_n)\to\cdots.
$$
Since $_{\Lambda}M$ is a generator (that is, $\add(_{\Lambda}\Lambda)\subseteq \add(M)$), the $A\opp$-module $M_A =\Hom_{\Lambda}(\Lambda, M)$ is projective. Therefore there is a series of isomorphisms of ${\Lambda}$-$A$-bimodules:
$$\begin{array}{ll} M\otimes_A D(A)& \simeq \Hom_{\Lambda}(\Lambda, M)\otimes_AD(\Hom_{\Lambda}(M,M))\\ & \simeq D\Hom_{A^{\opp}}\big(\Hom_{\Lambda}(\Lambda,M),\Hom_{\Lambda}(M,M)\big) \\ &\simeq D\Hom_{\Lambda}(M, \Lambda)\\ & \simeq \nu_{\Lambda}(M).\end{array}$$
Here, the second isomorphism follows from \cite[Proposition 20.11, p.243 ]{AF} and the third one is referred to \cite[Lemma 2.2(2)]{xiartm} for hints. Thus the sequence $(\ast\ast)$ is isomorphic to the following sequence
$$0\lra\Hom_{\Lambda}(\nu_{\Lambda}(M), M)\lra \Hom_{\Lambda}(\nu_{\Lambda}(M), I_0)\lra\cdots\lra \Hom_{\Lambda}( \nu_{\Lambda}(M), I_n)\lra\cdots.$$
Consequently, $\Ext_A^n(D(A), A)\simeq \Ext_{\Lambda}^n(\nu_{\Lambda}(M), M)$ as $A$-modules for all $n\geq 0$. Thus $\Ext_A^n(D(A), A)=0$ if and only if $\Ext_{\Lambda}^n(\nu_{\Lambda}(M), M)=0$. $\square$

\medskip
{\bf Proof of Theorem \ref{Main result}.}  Let $A$ be an Artin algebra, $M$ a generator-cogenerator for $A\modcat$, and  $B:=\End_A(M)$. Suppose $\Ext_A^n(M\oplus \nu_A(M), M)=0$ for all $n\geq 1$. By Lemma \ref{WG}, $\dm(B)=\infty$ and $\Ext_B^n(D(B), B)=0$ for all $n\geq 1$. Suppose that the algebra $B$ is virtually Gorenstein. Then $B$ is self-injective by Theorem \ref{DM+EXT}.  Since  ${_A}M$ is also a generator for $A\Modcat$,  the functor $\Hom_A(M,-): A\Modcat\to B\Modcat$ is fully faithful. Let $f: M\to I_0$ be an injective envelope of $_AM$. Since ${_A}M$ is a generator-cogenerator, the $B$-module $\Hom_A(M, I_0)$ is projective-injective and $\Hom_A(M, f)$ is an injective envelope of $_BB$. Thus $\Hom_A(M,M)={}_BB\simeq {}_B\Hom_A(M, I_0)$. This implies $M\simeq I_0$ as $A$-modules.  In particular, ${_A}M$ is injective. Since ${_A}M$ is a generator, it follows from $\add(_AA)\subseteq \add(_AM)$ that the algebra $A$ itself is self-injective, and therefore ${_A}M$ is also projective.  $\square$

\medskip
Theorem \ref{Main result} involves both infinite dominant dimensions and orthogonality of modules. We introduce the following property $(\diamondsuit)$ for an algebra $A$ and show that this property is preserved by taking tensor products of algebras over a field.

$(\diamondsuit)$:\;\; $\dm(A)=\infty$ and $\Ext_A^n(D(A), A)=0$ for all $n\geq 1$.

\begin{Prop}\label{TAG}
Let $A$ and $B$ be finite-dimensional algebras over a field $k$ and let $C:=A\otimes_kB$ be the tensor product of $A$ and $B$ over $k$. Then $A$ and $B$ satisfy the property $(\diamondsuit)$ if and only if so does $C$.
\end{Prop}

{\it Proof.} When either $A$ or $B$ is zero, Proposition \ref{TAG} holds trivially. So, we assume that
both $A$ and $B$ are nonzero. By \cite[Lemma 3]{Muller}, $\dm(C)=\min\{\dm(A), \dm(B)\}$. This implies that $\dm(C)=\infty$ if and only if $\dm(A)=\infty=\dm(B)$. Note that $D(C)\simeq D(A)\otimes_kD(B)$ as $C$-$C$-bimodules. Since $k$ is a field, it follows from \cite[Chapter XI, Theorem 3.1]{CE} that, for all $n\in\mathbb{N}$, there are isomorphisms of $k$-modules: $$(\sharp)\quad \Ext_C^n(D(C),C)\simeq\Ext_C^n(D(A)\otimes_kD(B), A\otimes_k B)\simeq
\bigoplus_{p, q\geq 0, p+q=n}\Ext_A^p(D(A), A)\otimes_k\Ext_B^q(D(B), B).$$
Suppose $\dm(A)=\infty=\dm(B)$. Then $A$ and $B$ have finite-dimensional, projective-injective, nonzero modules, and therefore $\Hom_A(D(A),A)\neq \Hom_B(D(B), B)$.
It follows from $(\sharp)$ that $\Ext_C^n(D(C), C)$ = $0$ for all $n\geq 1$ if and only if
$\Ext_A^m(D(A), A)=0=\Ext_B^m(D(B), B)$ for all $m\geq 1$. Thus $A$ and $B$ satisfy the property $(\diamondsuit)$ if and only if so does $C$. $\square$

\medskip
Next, we give a generalization of Theorem \ref{DM+EXT} in the case of finite-dimensional algebras.

\begin{Theo}
Suppose that $A$ is a finite-dimensional algebra over a field $k$ with $\dm(A)=\infty$ and $\Ext_A^n(D(A), A)=0$ for all $n\geq 1$. If $A$ is isomorphic to the tensor product of virtually Gorenstein algebras, then $A$ is self-injective.
\end{Theo}

{\it Proof.} Suppose  $A\simeq A_1\otimes_k \cdots \otimes_k A_m$, where $A_i$ is a virtually Gorenstein algebra for $1\leq i\leq m$. Since $\dm(A)=\infty$ and $\Ext_A^n(D(A), A)=0$ for all $n\geq 1$, we see from Proposition \ref{TAG} that, for $1\leq i\leq m$, $\dm(A_i)=\infty$ and $\Ext_{A_i}^j(D(A_i), A_i)=0$ for all $j\geq 1$. By Theorem \ref{DM+EXT}, $A_i$ is self-injective.
Note that the tensor product of finitely many, self-injective algebras is again a self-injective algebra. Thus $A$ is self-injective. $\square$

\medskip
Finally, we give a way to get virtually Gorenstein algebras.

Let $B\subseteq A$ be an extension of algebras, that is, $B$ is a subalgebra of the algebra $A$ with the same identity. An extension $B\subseteq A$ is called a \emph{Frobenius} extension if ${_B}A$ is a finitely generated projective $B$-module and ${_A}A_{B}\simeq \Hom_B(A,B)$ as $A$-$B$-bimodules (\cite{kasch}). This is equivalent to saying that  $A{_B}$ is a finitely generated projective $B\opp$-module and $_{B}A{_A}\simeq \Hom_{B\opp}(A,B)$ as $B$-$A$-bimodules. Given a Frobenius extension $B\subseteq A$, it is known that the restriction functor ${_B}(-): A\Modcat\to B\Modcat$ and the induction function $A\otimes_B-: B\Modcat\to A\Modcat$ are mutually adjoint, and thus preserve projective (respectively, injective) modules. However, they do not detect projective (respectively, injective) modules in general. For example, the inclusion $k\subseteq k[x]/(x^n)$ with $k$ a field and $n\ge 2$ is a Frobenius extension, the restriction of every $k[x]/(x^n)$-module is a projective $k$-module, but the module itself may not be a projective $k[x]/(x^n)$-module. So, to establish close relation between modules over $A$ and $B$, we focus on two classes of special Frobenius extensions. For more examples of Frobenius extensions, we refer to \cite{Kad}.

An extension $B\subseteq A$ of algebras is called a \emph{separable} extension if the multiplication $A\otimes_BA\to A$ is a split surjection of $A$-$A$-bimodules; a \emph{semisimple} extension if the multiplication map $A\otimes_BX\to X$ is split surjective for any $A$-module $X$; and a \emph{split} extension if the inclusion $B\to A$ is a split injection of $B$-$B$-bimodules. Clearly, separable extensions are semisimple. For a semisimple extension $B\subseteq A$, if $X$ is an $A$-module, then $X$ is isomorphic to a direct summand of the $A$-module $A\otimes_BX$; for a split extension, if $Y$ is a $B$-module, then $Y$ is isomorphic to a direct summand of the $B$-module ${_B}A\otimes_BY$.

\begin{Prop}\label{frob-ext}
Let  $B\subseteq A$ be a Frobenius extension of algebras.

$(1)$ If the extension is semisimple and $B$ is virtually Gorenstein, then $A$ is virtually Gorenstein.

$(2)$ If the extension is split and $A$ is virtually Gorenstein, then $B$ is virtually Gorenstein.
\end{Prop}

{\it Proof.} Let $A\GP$ denote the category of all Gorenstein-projective $A$-modules. Define
$$F={_B}(-): A\Modcat\lra B\Modcat,\quad G=A\otimes_B-: B\Modcat\lra A\Modcat.$$
Clearly, $(G,F)$ is an adjoint pair. Since $B\subseteq A$ is a Frobenius extension, $G$ is naturally isomorphic to the coinduction functor $\Hom_B(A,-)$. This implies that $(F, G)$ is also an adjoint pair. It is not difficult to see that $F$ and $G$ can be restricted to mutually adjoint functors between $\GP{A}$ and $\GP{B}$. 
Since $F$ and $G$ are exact and preserve projective modules, they automatically induce mutually adjoint triangle functors (still denoted by $F$ and $G$) between the stable category of $\GP{A}$ and the one of $\GP{B}$:
$$
F: A\mbox{-}\underline{\rm GProj}\lra B\mbox{-}\underline{\rm GProj}\,, \quad G: B\mbox{-}\underline{\rm GProj}\lra A\mbox{-}\underline{\rm GProj}.
$$
Note that, for an Artin algebra $\Lambda$, the category $\Lambda\mbox{-}\underline{\rm GProj}$ is a compactly generated triangulated category (for example, see \cite[Theorem 6.6]{Bel}). Now, we denote by $A\mbox{-}\underline{\rm GProj}^{\rm c}$ and $B\mbox{-}\underline{\rm GProj}^{\rm c}$ the full subcategories of $A\mbox{-}\underline{\rm GProj}$ and $B\mbox{-}\underline{\rm GProj}$ consisting of all compact objects, respectively. Recall that an object of a triangulated category $\mathscr{C}$ with coproducts (indexed by sets) is said to be \emph{compact} if the functor $\Hom_\mathscr{C}(X,-)$ from $\mathscr{C}$ to the category of abelian groups commutes with coproducts. For the convenience of the reader, we mention two general results:

$(a)$ Given an adjoint pair $(L, R)$ of triangle functors $L:\mathscr{C}\to\mathscr{D}$ and $R:\mathscr{D}\to\mathscr{C}$ between triangulated categories $\mathscr{C}$ and $\mathscr{D}$ with coproducts, if $R$ commutes with coproducts, then $L$ preserves compact objects. This fact is easy to see by definition.

$(b)$ Let $\Lambda\mbox{-}\underline{\rm Gproj}$ be the full subcategory of $\Lambda\mbox{-}\underline{\rm GProj}$ consisting of modules isomorphic to finitely generated Gorenstein-projective $\Lambda$-modules. Then $\Lambda\mbox{-}\underline{\rm Gproj}\subseteq \Lambda\mbox{-}\underline{\rm GProj}^{\rm c}$,
and the equality holds if and only if $\Lambda$ is virtually Gorenstein (see \cite[Theorem 8.2 (i) and (iv)]{Bel}).

Since $F$ and $G$ are mutually adjoint and commute with coproducts, they can be restricted to
triangle functors between $A\mbox{-}\underline{\rm GProj}^{\rm c}$ and $B\mbox{-}\underline{\rm GProj}^{\rm c}$ by $(a)$. Thus we obtain mutually adjoint pairs:
$$F: A\mbox{-}\underline{\rm GProj}^{\rm c}\lra B\mbox{-}\underline{\rm GProj}^{\rm c}\,, \quad G: B\mbox{-}\underline{\rm GProj}^{\rm c}\lra A\mbox{-}\underline{\rm GProj}^{\rm c}.
$$

$(1)$ Suppose that the extension $B\subseteq A$ is separable and $B$ is virtually Gorenstein.
Let $X\in A\mbox{-}\underline{\rm GProj}^{\rm c}$. Then $F(X)\in B\mbox{-}\underline{\rm GProj}^{\rm c}$ and $GF(X)\in A\mbox{-}\underline{\rm GProj}^{\rm c}$. Since $B$ is virtually Gorenstein, $F(X)\in B\mbox{-}\underline{\rm Gproj}$. It follows that $GF(X)\in A\mbox{-}\underline{\rm Gproj}$. Since the extension $B\subseteq A$ is semisimple, $X$ is isomorphic to a direct summand of $GF(X)$. As $A\mbox{-}\underline{\rm Gproj}$ is closed under direct summands in $A\mbox{-}\underline{\rm GProj}$, we have $X\in A\mbox{-}\underline{\rm Gproj}$. Thus $A\mbox{-}\underline{\rm Gproj} = A\mbox{-}\underline{\rm GProj}^{\rm c}$. It follows from $(b)$ that $A$ is virtually Gorenstein.

$(2)$ This can similarly be shown by applying $(b)$ and the fact that, for a split extension $B\subseteq A$, each $B$-module $Y$ is isomorphic to a direct summand of the $B$-module $G(Y)$. $\square$

\medskip
{\bf Acknowledgement.}
The research work was supported partially by the National Natural Science Foundation of China (Grant 12031014).

{\footnotesize
\smallskip
Hongxing Chen,

School of Mathematical Sciences  \&  Academy for Multidisciplinary Studies, Capital Normal University, 100048
P. R. China

{\tt Email: chenhx@cnu.edu.cn}

\smallskip
Changchang Xi,

School of Mathematical Sciences, Capital Normal University, 100048
Beijing, P. R. China

{\tt Email: xicc@cnu.edu.cn}
}


{\footnotesize
\begin{thebibliography}{99}
\bibitem{AF} {{\sc F. W. Anderson} and {\sc K. R. Fuller}, \emph{Rings and categories of modules}, Second Edition.
Graduate Texts in Mathematics \textbf{13}, Springer-Verlag 1992.}

\bibitem{AR2} {{\sc M. Auslander} and {\sc I. Reiten}, Applications of contravariantly finite subcategories,
{\it Adv. Math.} {\bf 86} (1991) 111-152.}

\bibitem{ars}{{\sc M. Auslander, I. Reiten} and {\sc S. Smal{\o}},
\emph{Representation theory of Artin algebras}, Cambridge Studies in Advanced Mathematics \textbf{36}.
Cambridge University Press, Cambridge, 1995.  }

\bibitem{Bel} {{\sc A. Beligiannis}, Cohen-Macaulay modules, (co)torsion pairs and virtually Gorenstein
algerbas, {\it J. Algebra} {\bf 288} (2005) 137-211.}

\bibitem{Bel2} {{\sc A. Beligiannis}, On algebras of finite Cohen-Macaulay type,
{\it Adv. Math.} {\bf 226} (2011) 1973-2019.}

\bibitem{BK} {{\sc A. Beligiannis} and {\sc H. Krause}, Thick subcategories and virtually Gorenstein algerbas,
{\it Illinois J. Math}. {\bf 52} (2008) 551-562.}

\bibitem{BI}{{\sc A. Beligiannis} and {\sc I. Reiten},
Homological and homotopical aspects of torsion theories, {\it Mem.
Amer. Math. Soc}. {\bf 188} (883) (2007) 1-207.}

\bibitem{CE} {{\sc H.Cartan} and {\sc S. Eilenberg}, Homological algebra. Princeton University Press,  Princeton, New Jersey, 1956.}


\bibitem{CFKKY} {{\sc H.X. Chen, M. Fang, O. Kerner, S. Koenig} and {\sc K. Yamagata}},
Rigidity dimension of algebras, \emph{Math. Camb. Phil. Soc.} \textbf{170} (2021) 417-443.

\bibitem{xcf}{{\sc H. X. Chen}, {\sc M. Fang} and {\sc C. C. Xi}, Tachikawa's second conjecture, derived
recollements, and gendo-symmetric algebras, \emph{Compos. Math.} \textbf{160} (11) (2024) 2704-2734.}

\bibitem{xc6}{{\sc H. X. Chen} and {\sc C. C. Xi}, Dominant dimensions, derived equivalences
and tilting modules, \emph{Isr. J. Math.} \textbf{215} (2016) 349-395.}

\bibitem{xc3}{{\sc H. X. Chen} and {\sc C. C. Xi}, Homological theory of self-orthogonal modules,
\emph{Trans. Amer. Math. Soc.} \textbf{378} (10) (2025) 7287-7335.  arXiv:2208.14712.}

\bibitem{FK11} {{\sc M. Fang} and {\sc S. Koenig}},
 Endomorphism algebras of generators over symmetric algebras,
 \textit{J. Algebra} \textbf{332} (2011) 428-433.

\bibitem{FK14} {{\sc M. Fang} and {\sc S. Koenig}},
 Gendo-symmetric algebras, canonical comultiplication, bar
 cocomplex and dominant dimension,
\textit{Trans. Amer. Math. Soc.} \textbf{368} (2016) 5037-5055.

\bibitem{OI} {\sc O. Iyama},
Higher-dimensional Auslander-Reiten theory on maximal orthogonal subcategories,
\textit{Adv. Math.} \textbf{210} (2007) 22-50.

\bibitem{Kad}{{\sc L. Kadison}, {\em New examples of Frobenius extensions}, University Lecture Series \textbf{14}, American Mathematical Society, 1999.}

\bibitem{kasch}{{F. Kasch,} Projektive Frobenius-Erweiterungen. \emph{S.-B. Heidelberger Akad Wiss Math-Nat Kl,}  \textbf{1960(61)} (1960/61) 87-109.}

\bibitem{KS} {{\sc H. Kruase} and {\sc S. {\O}yvind}},
Applications of cotorsion pairs,
\textit{J. London Math. Soc. (2)} \textbf{68} (2003) 631-650.

\bibitem{mar}{{\sc R.  Marczinzik}, Upper bounds for the dominant dimension of Nakayama and related algebras, \emph{J. Algebra} \textbf{496} (2018) 216-241.}

\bibitem{Mo58a}{{\sc K. Morita}, Duality for modules and its applications of the theory of rings with minimum condition,
\emph{Sci. Rep. Tokyo Kyoiku Daigaku Sec. A} \textbf{6} (1958) 83-142.}

\bibitem{Muller} {{\sc B. M\"{u}ller}, The classification of algebras by dominant dimension, \emph{Canad. J. Math.} \textbf{20} (1968) 398-409.}

\bibitem{Nakayama}{{\sc T. Nakayama}, On algebras with complete homology, \emph{Abh. Math. Sem. Univ.
Hamburg} \textbf{22} (1958) 300-307.}

\bibitem{nrtz}{{\sc V. C. Nguyen, I. Reiten, G. Todorov,} and {\sc S. J. Zhu}, Dominant dimension and tilting modules,
\emph{Math. Z.} \textbf{292}(3-4) (2019) 947-973.}

\bibitem{Tac} {{\sc H. Tachikawa}, \emph{Quasi-Frobenius rings and generalizations}. Notes by Claus Michael Ringel.
Lecture Notes in Mathematics \textbf{351}, Berlin-Heidelberg New York 1973.}

\bibitem{xiartm}{{\sc C. C. Xi}, The relative Auslander-Reiten theory of modules. Preprint is
available at: https://www.wemath.cn/
$^\sim$ccxi/Papers/Articles/rart.pdf.}

\bibitem{xi-2}{{\sc C. C. Xi}, Frobenius bimodules and flat-dominant dimensions,
\emph{Sci. China Math.} \textbf{64} (2021) 33-44.}

\end{thebibliography}}
\end{document}